\documentclass[11pt]{article}

\usepackage{amsfonts} 

\newcommand{\C}{\mathbb{C}}
\newcommand{\F}{\mathbb{F}}
\newcommand{\R}{\mathbb{R}}

\newcommand{\Z}{\mathbb{Z}}

\newcommand{\TT}{\mathbb{T}}

\newtheorem{thm}{Theorem}
\newtheorem{defn}[thm]{Definition}

\newtheorem{coro}[thm]{Corollary}

\newtheorem{rem}[thm]{Remark}

\newcommand{\dd}{\;\mathrm{d}}
\newcommand{\e}{\mathrm{e}}
\newcommand{\binom}[2]{{#1\choose#2}} 
\newcommand{\pf}{\noindent {\bf PROOF.} \quad}
\newcommand{\qed}{$\Box$}
\newcommand{\Vol}{\mathrm{Vol}}

\newcommand{\tr}{\mathrm{tr}} 
\newcommand{\re}{\mathop{\mathrm{Re}}} 
\newcommand{\im}{\mathop{\mathrm{Im}}} 
\newcommand{\ii}{\mathrm{i}} 
\begin{document}
\title{Mahler measure under variations of the base group}
\author{Oliver T. Dasbach
{\footnote{The first author was supported in part by NSF grants DMS-0306774 and DMS-0456275 (FRG)}}\\ 
\small {Louisiana State University}\\
\small {Department of Mathematics}\\
\small {Baton Rouge, LA 70803}\\
\small {e-mail:  \texttt{kasten@math.lsu.edu}} 
\and
Matilde N. Lalin
{\footnote{The second author was partially supported by the NSF under DMS-0111298.}}\\
\small{University of British Columbia}\\
\small{Department of Mathematics}\\
\small{Vancouver, BC V6T 1Z2, Canada}\\
\small{e-mail: \texttt{mlalin@math.ubc.ca}}}

\maketitle

\begin{abstract}
We study properties of a generalization of the Mahler measure to elements in group rings, in terms of the L\"uck-Fuglede-Kadison determinant. Our main focus is the variation of the Mahler measure when the base group is changed. In particular, we study how to obtain the Mahler measure over an infinite group as limit of Mahler measures over finite groups, for example, in the classical case of the free abelian group or the infinite dihedral group, and others.  
\end{abstract}



\section{Introduction}
The Mahler measure of a Laurent polynomial $P \in \C[x_1^{\pm1}, \dots, x_n^{\pm1}]$,  is defined by
\begin{eqnarray}
m(P)& := &\frac{1}{(2  \pi \ii)^n} \int_{\TT^n} \log |P(x_1, \dots, x_n) | \frac{ \dd x_1}{x_1} \dots \frac{ \dd x_n}{x_n} \\
& = & \int_0^1 \dots \int_0^1 \log |P(\e^{2 \pi \ii t_1}, \dots , \e^{2 \pi \ii t_n}) | \dd t_1 \dots \dd t_n. 
\end{eqnarray}

In particular, Jensen's formula implies, for a one-variable polynomial that 
\begin{equation}
m(P)= \log|a| +\sum_i \log \max \{|\alpha_i|, 1 \} \qquad \mathrm{for} \quad P(x)=a \prod_i (x-\alpha_i).
\end{equation}

In this article we use some ideas of Rodriguez-Villegas \cite{RV:ModularMahler} to generalize the definition of the Mahler measure to certain elements of group rings. Our construction yields the logarithm of the L\"uck-Fuglede-Kadison determinant as defined by L\"uck in \cite{LueckBook:L2invariants}. 

We relate the Mahler measure for elements in the group ring of a finite group to the characteristic polynomial of the adjacency matrix of a (weighted) Cayley graph and also to characters of the group (Theorem \ref{thm:finitegroup}). We obtain certain
rationality results for the Mahler measure of group ring elements of the form $1-\lambda P$. In Section \ref{sect:finite_approximation} we  recover a result of L\"uck \cite{LueckBook:L2invariants} that the (classical) Mahler measure is limit of Mahler measures over finite abelian groups. 
We extend this result to a more general case. In particular, we complement it in Theorem \ref{thm:generalapprox} and Corollary \ref{everything} by proving an analogous property for the infinite dihedral group, $PSL_2(\Z)$, and other infinite groups. Finally, we study how the Mahler measure changes when we vary the base group. In particular we focus our attention to the dihedral group $D_m$ versus the abelianization $\Z/m\Z \times \Z/2\Z$ in the general case (Theorem \ref{thm:dihedral}).

\section{The general technique}
In this section we will review the techniques and one example of Rodriguez-Villegas \cite{RV:ModularMahler}. The section after that will point out the connections to the combinatorial $L^2$-torsion of L\"uck \cite{LueckBook:L2invariants}. Let $P \in \C[x^{\pm 1},y^{\pm 1}]$  be a polynomial, $\lambda \in \C$, and let $P_\lambda(x,y)= 1 - \lambda P(x,y)$. Assume that $P$ is a reciprocal polynomial i.e., $P(x,y) = \overline{P\left(x^{-1},y^{-1}\right)}$. For instance, as in \cite{RV:ModularMahler}, we may think of  $P(x,y) = x+x^{-1}+y+y^{-1}$. 

We define $m(P,\lambda) :=m(P_\lambda)$. Then

\begin{eqnarray}
m(P,\lambda) &=& \frac{1}{(2\pi\ii)^2} \int_{\TT^2} \log|1-\lambda P(x,y)| \frac{\dd x}{x}\frac{\dd y }{y}.\\
&=& \int_0^1 \int_0^1 \log |P(\e^{2 \pi i t}, \e^{2 \pi i s})| \, {\dd t} \, {\dd s}.
\end{eqnarray} 
 
Since $P(x,y)$ is a continuous function defined in the torus $\TT^2$, it is bounded and  $|\lambda P(x,y)|<1$ for $\lambda$ small enough and $x,y \in \TT^2$. For instance, it is easy to see that we need to take $|\lambda|<\frac{1}{4}$ in the example. If $\lambda$ is real, $1- \lambda P(x,y)>0$ on the torus $\TT^2$.


Let
\begin{eqnarray}
u(P,\lambda) &=& \frac{1}{(2\pi\ii)^2} \int_{\TT^2} \frac{1}{1-\lambda P(x,y)} \frac{\dd x}{x}\frac{\dd y }{y} \nonumber \\
&=& \sum_{n=0}^\infty \lambda^n \frac{1}{(2\pi\ii)^2} \int_{\TT^2} P(x,y)^n \frac{\dd x}{x} \frac{\dd y}{y}= \sum_{n=0}^\infty a_n \lambda^n,
\end{eqnarray}
where
\[\frac{1}{(2\pi\ii)^2} \int_{\TT^2} P(x,y)^n \frac{\dd x}{x} \frac{\dd y}{y}= \left[P(x,y)^n\right]_0 = a_n\]
corresponds to the constant coefficient of the $n$-th power of $P$. Now we can write
\begin{eqnarray} \label{eq:defn}
 \tilde{m}(P,\lambda) &=& \frac{1}{(2\pi\ii)^2} \int_{\TT^2} \log(1-\lambda P(x,y)) \frac{\dd x}{x}\frac{\dd y }{y}\nonumber \\ 
\label{MMeasureSeries} 
&=&- \int_0^\lambda (u(t)-1) \frac{\dd t}{t}= - \sum_{n=1}^\infty \frac{a_n \lambda^n}{n}.
\end{eqnarray}
In general,
\begin{equation}
m(P,\lambda) = \re( \tilde{m}(P,\lambda)).
\end{equation}

In our particular example, $\tilde{m}(\lambda)$ is real. For $\lambda$ real, $a_n = 0$ for $n$ odd and 
\[ a_{2m} = \sum_{j=0}^m \frac{(2m)!}{j!j!(m-j)!(m-j)!} = \binom{2m}{m}\sum_{j=0}^m \binom{m}{j}^2 = \binom{2m}{m}^2.\]

We observe that the right term in equation (\ref{eq:defn}) makes sense in more general contexts. More precisely, we can make the following:
\begin{defn} \label{def:MahlerM}
Let $ \F_{x_1,\dots,x_l}$ be the free group in $x_1, \dots, x_l$ and let $N$ be a normal subgroup (of relations). 
Let $\Gamma = \F_{x_1,\dots, x_l}/N$ be the quotient group and consider the group ring $\C \Gamma$. 

Given a 
\[Q= Q(x_1, \dots x_l)=\sum_{g \in \Gamma} c_g g \in \C\Gamma,\] its reciprocal is 
\[Q^* = \sum_{g \in \Gamma} \overline{c_g} g^{-1} \in \C \Gamma.\] 
  
Let $P=P(x_1, \dots , x_l) \in \C\Gamma$ reciprocal,   
(i.e., $P=P^*$). Let $|\lambda| < \frac{1}{k}$ where $k$ is the $l_1$-norm of the coefficients of 
$P$ as an element in $\C \Gamma$. Then we define the Mahler measure of $P_\lambda(x_1,\dots, x_l) = 1-\lambda P(x_1,\dots x_l)$ by 
\begin{equation}
m_{\Gamma}(P, \lambda) = - \sum_{n=1}^\infty \frac{a_n \lambda^n}{n},
\end{equation}
where $a_n$ is the constant coefficient of the $n$th power of $P(x_1,\dots, x_l)$, in other words,
\begin{equation}
a_n = \left[P(x_1,\dots, x_l)^n\right]_0. 
\end{equation}

We will also use the notation
\[u_{\Gamma}(P, \lambda)=\sum_{n=0}^{\infty} a_n \lambda^n\]
for the generating function of the $a_n$.
\end{defn}

We will confine ourselves mainly to the two variable case and 
we will keep the notation $m(P,\lambda)$ for $m_{\Z\times\Z}(P,\lambda)$.

\begin{rem}
Let us observe that we can extend this definition for any polynomial $Q(x_1,\dots,x_l) \in \C\Gamma$. It suffices to write

\[QQ^*=\frac{1}{\lambda}\left( 1-\left(1-\lambda QQ^*\right)\right)\]
for $\lambda$ real and positive and $1/\lambda$ larger than the length of $QQ^*$. 
Define
\begin{equation}
m_\Gamma(Q) = -\frac{\log \lambda}{2} - \sum_{n=1}^\infty \frac{b_n}{2n},
\end{equation}
where
\[ b_n = \left[(1-\lambda QQ^*)^n\right]_0.\]

This sum might not necessarily converge. When it is convergent then
it is clear that $m_\Gamma(Q)$ is well-defined as it does not depend on the choice of $\lambda$: write
\[ (1-(1-\lambda \mu x))= \lambda \mu x = (1- (1- \lambda x)) \mu\]
then
\[ \log(1-(1-\lambda \mu x)) - \log(\lambda \mu) = \log(1-(1-\lambda  x)) + \log \mu -\log(\lambda \mu)\]
\[ = \log(1-(1-\lambda  x) )- \log \lambda.\]
\end{rem}

\section{L\"uck's combinatorial $L^2$-torsion.}

The approach of Rodriguez-Villegas to the Mahler measure, as explained above, has a powerful generalization
in the theory of combinatorial $L^2$-torsion of L\"uck \cite{LueckBook:L2invariants}.

For example, let $K$ be a knot and \[\Gamma=\pi_1(S^3 \setminus K)=\langle x_1, \dots, x_g \,|\, r_1, \dots, r_{g-1} \rangle\]
be a presentation of the fundamental group of the knot complement, such that no $x_j$ represents a trivial element in $\Gamma$.
For a square matrix $M$ with entries in $\C \Gamma$ the trace $\tr_{\C \Gamma}(M)$
denotes the coefficient of the unit element in the sum of the diagonal elements. 
The matrix $A^*=(a^*_{j,i})$ is the transpose of $A=(a_{i,j})$ with 
\[\left ( \sum_{g \in \Gamma} c_g g \right ) ^* = \sum_{g \in \Gamma} \overline {c_g} g^{-1}.\]

Let \[F=\left ( 
\begin{array}[pos]{ccc}
\frac{\partial r_1} {\partial x_1} & \dots & \frac{\partial r_1} {\partial x_g}\\
\vdots &\ddots & \vdots\\
\frac{\partial r_{g-1}} {\partial x_1} & \dots & \frac{\partial r_{g-1}} {\partial x_g}\\
\end{array}
\right )\]
be the Fox matrix (e.g. \cite{BZ}) of the presentation.
The matrix $F$ is a $(g-1)\times g$ matrix with entries in the group ring $\C \Gamma$.
We obtain a $(g-1) \times (g-1)$-matrix $A$ by deleting one of the columns of $F$.

\begin{thm}[L\"uck \cite{LueckBook:L2invariants}] \label{LuecksTheorem} Suppose $K$ is a hyperbolic knot.
Then, for $k$ sufficiently large

\begin{equation} \label{LuecksFormula}
\frac 1 {3 \pi} \Vol(S^3 \setminus K)= 2 (g-1) \ln (k) - \sum_{n=1}^{\infty} \frac 1 n \tr_{\C \Gamma} \left ( ( 1-k^{-2} A A^*)^n \right). 
\end{equation}
\end{thm}

In the special case that $A$ has entries in $\C [t,t^{-1}]$ it follows from (\ref{MMeasureSeries}) that
the right-hand side of (\ref{LuecksFormula}) is
$2 m(\det(A))$.

\section{Spectra of weighted Cayley graphs}

It will turn out that the Mahler measure is intimately related to the spectral theory of the
Cayley graph of the base group relative to the generating set given by the terms of the polynomial.

We recall a theorem of Babai \cite{Babai:SpectraCayleyGraphs}. 

For a finite group $\Gamma$ of order $m$ let $\alpha$ be a function $\alpha: \Gamma \rightarrow \C$ such that 
$\alpha(g)=\overline{\alpha(g^{-1})}$ for all elements $g \in \Gamma$. The weighted Cayley graph of $\Gamma$ with 
respect to $\alpha$ is a graph with vertices $g_1, \dots, g_m$. Two vertices $g_i$ and $g_j$ are
connected by a (directed) edge of weight $\alpha(g_i^{-1} g_j)$. Weight $0$ means that there is no edge between
$g_i$ and $g_j$. The weighted adjacent matrix $A(\Gamma, \alpha)$ of the weighted Cayley graph of $\Gamma$ with respect to
$\alpha$ is the matrix $A$ with entries $a_{i,j}=\alpha(g_i^{-1} g_j)$. Thus the matrix $A$ is hermitian: $A=A^*$.

The spectrum of the Cayley graph of $\Gamma$ with respect to $\alpha$ is the spectrum of 
$A(\Gamma,\alpha)$.

The spectrum can be computed from the irreducible characters of $\Gamma$ \cite{Babai:SpectraCayleyGraphs}.
Let $\chi_1, \dots, \chi_h$ be the irreducible characters of $\Gamma$ of degrees $n_1, \dots, n_h$.

\begin{thm}[Babai \cite{Babai:SpectraCayleyGraphs}] \label{TheoremBabai}
The spectrum of $A(\Gamma,\alpha)$ can be arranged as
\[{\cal S}=\left \{\sigma_{i,j}: i=1,\dots,h; j=1, \dots,n_i \right \},\]
such that 
$\sigma_{i,j}$ has multiplicity $n_i$ and
\begin{equation}
\sigma_{i,1}^t + \dots + \sigma_{i,n_i}^t = \sum_{g_1, \dots, g_t \in \Gamma} \left ( \prod_{s=1}^t \alpha(g_s) \right ) \chi_i \left (\prod_{s=1}^t g_s \right ).
\end{equation}

\end{thm}

\begin{rem}
Notice that Theorem \ref{TheoremBabai} allows us to compute the spectrum of $A$ from the irreducible
characters of the base group. 

Indeed, by Theorem \ref{TheoremBabai},
\begin{eqnarray}
\tr(A^n) &=& \sum_{i=1}^h n_i \sum_{g_1, \dots, g_n \in \Gamma} \left ( \prod_{s=1}^n \alpha(g_s)\right) \chi_i\left(\prod_{s=1}^n g_s \right)\\
&=&  \sum_{i=1}^h n_i P^n\left[ \chi_i(g_1), \dots, \chi_i(g_m) \right].
\end{eqnarray}
The last term means that we first compute the $n$-th power of the element $P$ as a polynomial and then we evaluate each of the monomials in the representation $\chi_i$. This is not the usual evaluation in the algebra generated by the $x_i$, therefore we use the brackets to distinguish it.
\end{rem}

\section{The Mahler measure over finite groups}

Suppose  $P\in \C\Gamma$ with 
\[P = \sum_i (\delta_i S_i+ \overline{\delta_i} S_i^{-1})+ \sum_j \eta_j T_j\]
for some monomials $S_i \neq S_i^{-1}$, $T_j=T_j^{-1}$, $\delta_i \in \C$, $\eta_j \in \R$, and 
$S_i, T_j \in \Gamma$, $\Gamma$ a group. 
For reasons of simplicity we assume that the monomials of $P$ actually generate the whole group $\Gamma$, rather than a subgroup.

The coefficients $\delta_i, \eta_j$ induce a weight on the Cayley graph of $\Gamma$ with respect to the generating set
$\{S_i, T_j\}$ of $\Gamma$.
We are interested in the computation of $a_n=\left[P^n\right]_0$.
This corresponds to counting (with weights) the number of closed paths of length $n$ in the Cayley graph of $\Gamma$ that start and end at the identity. The weight of a closed path is the product over the weights of the
edges of the path. Therefore, the value of $a_n$ can be computed via the trace of the $n$-th power of $A$.

\begin{thm}\label{thm:finitegroup}
If the group $\Gamma$ is finite then 
\begin{equation}
m_{\Gamma} (P, \lambda)= \frac 1 {|\Gamma|} \log \det (I- \lambda A),
\end{equation}
where $A$ is the adjacency matrix of the Cayley graph (with weights) and $\frac 1 {\lambda}$ greater than the largest eigenvalue of the hermitian matrix $A$.

In particular $\exp(|\Gamma| m_{\Gamma}(P,\lambda))$ is a polynomial in $\lambda$. 
\end{thm}
This Theorem provides an analytic continuation of $m_{\Gamma}(P,\lambda)$ as a function of $\lambda$ in the complex plane $\C$ minus the eigenvalues of $A$.

\medskip

\pf  
Since the Cayley graph is vertex-transitive, this implies that
\[a_n = \frac 1 {|\Gamma|} \tr (A^n),\]
where  $\tr$ is the usual trace on matrices.

Now the Mahler measure is
\begin{eqnarray*} 
m_{\Gamma} (P, \lambda) &=&  - \sum_{n=1}^{\infty} \frac {a_n \lambda^n} n\\
&=& - \sum_{n=1}^{\infty} \frac 1 {|\Gamma|}\frac {\tr( (\lambda A)^n)} n\\ 
&=& \frac 1 {|\Gamma|} \log \det (I- \lambda A).
\end{eqnarray*}\qed


\begin{rem}
In the case of an arbitrary polynomial $Q$, 
\begin{equation}
m_{\Gamma}(Q)= \frac {1}{2 |\Gamma|} \log \det B, 
\end{equation}
where $B$ is the adjacency matrix corresponding to $Q Q^*$.
The Mahler measure $m_\Gamma(Q)$ is defined if  and only if $B$ is nonsingular.
\end{rem}

\begin{rem}
Recall that the generating function for the $a_n$ is
\[u_{\Gamma}(P,\lambda)=\sum_{j=0}^{\infty} a_n \lambda^n\]
and satisfies:
\[-\lambda \frac {\dd}{\dd \lambda} m_{\Gamma}(P,\lambda) = u_{\Gamma}(P,\lambda)-1.\]

Let $\sigma_i$ be the eigenvalues of $A$. By Theorem \ref{thm:finitegroup} we have:
\begin{equation}\label{eq:urational}
u_{\Gamma}(P,\lambda) = 1 - \frac {\lambda}{|\Gamma|} \frac {\frac {\dd}{\dd \lambda} \det(I - \lambda A)}{\det(I-\lambda A)} = \frac{1}{|\Gamma|} \sum_i \frac{1}{1-\lambda \sigma_i}.
\end{equation}

In particular, $u_{\Gamma}(P,\lambda)$ is a rational function in $\lambda$.
\end{rem}

\subsection{An example}

We will consider the example of the group 
\[\Z/3 \Z \times \Z/2 \Z = \langle x, y \,|\, x^3, y^2, [x,y]\rangle.\]
Let us compute the Mahler measure of the general polynomial
\[P_\lambda(x,y)= 1-\lambda \left(a +  b x+ \overline{b} x^{-1}+  c y  +  d y x + \overline{d} y x^{-1} \right).\]
Numbering the vertices $e, x, x^{-1}, y, yx, yx^{-1}$ as $1,\dots,6$, we get an adjacency matrix:
\[ \left(\begin{array}{cccccc}a & b & \overline{b} & c & d &\overline{d}\\\overline{b} & a & b & \overline{d} & c &d\\b& \overline{b}& a & d &\overline{d} &c\\c &d & \overline{d} & a& b &\overline{b}\\\overline{d} & c& d &\overline{b}& a &b\\d & \overline{d} & c & b & \overline{b} &a\end{array} \right) .\]
We see that the characteristic polynomial $\det(t I - A)$ is 
\[((t-a+c+\re(b-d))^2-3(\im(b-d))^2)((t-a-c+\re(b+d))^2-3(\im(b+d))^2)\]
\[\cdot(t-a-2\re(b)+c+2\re(d))(t-a-2\re(b)-c-2\re(d)).\]
Thus,
\[ m_{\Z/3\Z\times\Z/2\Z}(P,\lambda)  =  \frac{1}{6}\log ((1-\lambda(a-c-\re(b-d))^2-3\lambda(\im(b-d))^2)\]
\[+ \frac{1}{6}\log ((1-\lambda(a+c-\re(b+d))^2-3\lambda(\im(b+d))^2) \]
\[+\frac{1}{6}\left(  \log\left(1-\lambda (a+2\re(b)-c-2\re(d)) \right)+ \log\left(1-\lambda(a+2\re(b)+c+2\re(d)) \right) \right).\]

For instance, if $Q=1+x+y$ then $QQ^*= 3+x+x^{-1}+2y+yx+yx^{-1}$. It is easy to see that $\det B =81$. Therefore,
\[m_{\Z/3\Z\times\Z/2\Z}(1+x+y)= \frac{1}{12} \log 81 = \frac{\log 3}{3}.\]

\section{Abelian Groups}\label{sec:abel}

Suppose $\Gamma$ is a finite abelian group 
\[\Gamma = \Z/m_1 \Z \times \dots \times \Z/m_l \Z\]
of order
\[| \Gamma |=m=\prod_{j=1}^l m_j.\]

We obtain
\begin{coro}\label{coro}
\begin{equation}
m_{\Gamma}(P, \lambda) =\frac{1}{|\Gamma|}\log\left(\prod_{j_1, \dots, j_l}\left( 1- \lambda P (\xi^{j_1}_{m_1},\dots,\xi_{m_l}^{j_l}) \right) \right)
\end{equation}
where $\xi_k$ is a primitive root of unity.
\end{coro}
\pf This Corollary is consequence from the form of the representations for a finite abelian group. The characters are given by
\[\chi_{j_1, \dots, j_l} (s_1,\dots,s_l) = \xi_{m_1}^{j_1 s_1}\dots \xi_{m_l}^{j_ls_l}.\]
\qed

For the following result we recall some notation of Boyd and Lawton. For an integral vector $m=(m_1, \dots, m_l)$ let \[q(m) = \min \left\{ H(s) \, \left| \, s=(s_1, \dots, s_l)\in \Z^l,\, \sum_{i=1}^l m_i s_i =0 \right. \right\},\] where $H(s)= \max_{1\leq i \leq l} |s_i|$.

Thus, for a reciprocal polynomial $P \in \C [\Z^l]$ in the Laurent polynomial ring of $l$ variables, we have the following approximation.

\begin{thm}\label{thm:approx} For sufficiently small $\lambda$,
\begin{eqnarray}
\lim_{q(m) \rightarrow \infty}  m_{\Z/m_1\Z \times \dots \times \Z/m_l\Z}(P, \lambda)&=& m_{\Z^l} (P, \lambda).
\end{eqnarray}
\end{thm}
\pf Following Rodriguez-Villegas technique
\[m_{\Z^l} (P, \lambda)=  - \int_0^1 \dots \int_0^1 \log (1 - \lambda P(\e^{2 \pi i t_1},\dots, \e^{2 \pi i t_l})) \dd t_1 \dots \dd t_l.\]
Now the function  $\log (1 - \lambda P(\e^{2 \pi i t_1},\dots, \e^{2 \pi i t_l}))$ is continuous and with no singularities in $\TT^l$. Then we can compute the Riemann integral by taking limits on evaluations on the roots of unity. \qed
 
This Theorem should be compared to the result (for the one variable case) by L\"uck (page 478 of \cite{LueckBook:L2invariants}).

\section{Dihedral groups} 
 
Let us now consider $\Gamma = D_m$ the dihedral group, i.e.,
\[D_m= \langle \rho, \sigma \, | \, \rho^m, \sigma^2, \sigma \rho \sigma \rho  \rangle.\]
First assume that $m$ is odd, $m=2p+1$. The character  table (see, for instance \cite{Serre3}) for this group is 

\bigskip
\begin{center}
\begin{tabular}{c|c|rc}
& $\rho^k$ & $\sigma \rho^k$ & \\
\hline
$\chi_j$ & $\xi_m^{jk}+\xi_m^{-jk}$ & $0$ & $j=1, \dots p$\\
$\psi_1$ & $1$ & $1$ & \\
$\psi_2$ & $1$ & $-1$ & \\
\end{tabular} 
\end{center}

\bigskip

Assume that $P$ is reciprocal. In computing the Mahler measure of $1-\lambda P$ we need to determine
\[\tr(A^n) = P^n[1,\dots, 1]+P^n[1,\dots,1,-1, \dots, -1]\]
\[ + 2 \sum_{j=1}^p P^n\left[2,\xi_m^{j}+\xi_m^{-j}, \dots, \xi_m^{j(m-1)}+\xi_m^{-j(m-1)}, 0, \dots, 0\right].\]
Note that $P$ can be thought as a polynomial in two variables $\rho$ and $\sigma$. For the first two terms we can write
$P^n(1,1)+P^n(1,-1)$, then we obtain
\[=P^n(1,1)+P^n(1,-1)\]
\[ + 2 \sum_{j=1}^p \left(P^n\left[1,\xi_m^{j}, \dots, \xi_m^{j(m-1)}, 0, \dots, 0\right] +P^n\left[1,\xi_m^{-j}, \dots, \xi_m^{-j(m-1)}, 0, \dots, 0\right]\right) \]
\[=P^n(1,1)+P^n(1,-1)\]
\[ + 2 \sum_{j=1}^{m-1} P^n\left[1,\xi_m^{j}, \dots, \xi_m^{j(m-1)}, 0, \dots, 0\right]\]
\[=P^n(1,1)+P^n(1,-1) +  \sum_{j=1}^{m-1} \left(P^n\left(\xi_m^{j}, 1\right)+ P^n\left(\xi_m^{j}, -1\right) \right)\]
\begin{equation}\label{eq:dihedral}
=\sum_{j=1}^{m} \left(P^n\left(\xi_m^{j}, 1\right)+ P^n\left(\xi_m^{j}, -1\right) \right).
\end{equation}
At this point we should be more specific about the meaning of this formula. Given $P$ we consider its $n$th-power, then we write it as a combinations of monomials $\rho^k$ and $\sigma \rho^k$, and after that we evaluate each of the variables $\rho$, $\sigma$.

When $m =2p$, the character table is 

\bigskip
\begin{center}
\begin{tabular}{c|c|cc}
& $\rho^k$ & $\sigma \rho^k$ & \\
\hline
$\chi_j$ & $\xi_m^{jk}+\xi_m^{-jk}$ & $0$ & $j=1, \dots p-1$\\
$\psi_1$ & $1$ & $1$ & \\
$\psi_2$ & $1$ & $-1$ & \\
$\psi_3$ & $(-1)^k$ & $(-1)^k$ & \\
$\psi_4$ & $(-1)^k$ & $(-1)^{k+1}$ & \\
\end{tabular}
\end{center}

\bigskip
\noindent and we obtain a similar expression as in equation 
(\ref{eq:dihedral}).

Hence, we have proved:
\begin{thm} Let $P \in \C[D_m]$ be reciprocal. Then
\begin{equation}
\tr(A^n)= \sum_{j=1}^{m} \left(P^n\left(\xi_m^{j}, 1\right)+ P^n\left(\xi_m^{j}, -1\right) \right),
\end{equation}
where $P^n$ is expressed as a sum of monomials $\rho^k$, $\sigma \rho^k$ before being evaluated.
\end{thm}

On the other hand, for $\Gamma= \Z/m\Z\times \Z/2\Z = \langle x, y \,|\, x^m, y^2, [x,y]\rangle $,
\[ \tr (A^n) = \sum_{j=1}^m \left(P\left(\xi_m^j,1\right)^n + P\left(\xi_m^j,-1\right)^n\right).\]

From now on, in order to compare elements in the group rings of $D_m$ and $\Z/m\Z\times \Z/2\Z$ we will rename $x=\rho$ and $y=\sigma$ in $D_m$.

We have
\begin{thm}\label{thm:dihedral} Let 
\[P= \sum_{k=0}^{m-1} \alpha_k x^k  + \sum_{k=0}^{m-1} \beta_k y x^k\]
with real coefficients and reciprocal in $\Z/m\Z\times \Z/2\Z$ (therefore it is also reciprocal in $D_m$).
Then
\begin{equation}
m_{\Z/m\Z\times \Z/2\Z}(P,\lambda) = m_{D_m}(P,\lambda).
\end{equation}
\end{thm}
\pf It is important to note that being reciprocal in $\Z/m\Z\times \Z/2\Z$ is not the same as being reciprocal in $D_m$.

Being reciprocal in $\Z/m\Z\times \Z/2\Z$ means that $\alpha_k=\overline{\alpha_{m-k}}$ and $\beta_k = \overline{\beta_{m-k}}$. On the other hand, being  reciprocal in $D_m$ means that
$\alpha_k=\overline{\alpha_{m-k}}$ and $\beta_k$ is real. 

Then condition that $\alpha_k, \beta_k$ are real and the reciprocity in $\Z/m\Z\times \Z/2\Z$ imply that and $\alpha_k=\alpha_{m-k}$, $\beta_k=\beta_{m-k}$ and reciprocity in $D_m$. 

We will prove that the powers of $P$ are the same in both groups rings by induction. Assume that $P^n$ looks the same in both group rings and that it satisfies the conditions:
\[ P^n = \sum_{k=0}^{m-1} a_k x^k  + \sum_{k=0}^{m-1} b_k y x^k\]
with $a_k= a_{m-k}$ and $b_k= b_{m-k}$.
Then, in $\Z/m\Z\times \Z/2\Z$,
\[ P P^n= \sum_{k_1, k_2} \left( \alpha_{k_1} a_{k_2}+\beta_{k_1} b_{k_2}\right) x^{k_1+k_2} +
 \sum_{k_1, k_2} \left( \alpha_{k_1} b_{k_2}+\beta_{k_1} a_{k_2}\right) y x^{k_1+k_2}.\]
In $D_m$,
\[ P P^n= \sum_{k_1, k_2} \left( \alpha_{k_1} a_{k_2}+\beta_{m-k_1} b_{k_2}\right) x^{k_1+k_2} +
 \sum_{k_1, k_2} \left( \alpha_{m-k_1} b_{k_2}+\beta_{k_1} a_{k_2}\right) y x^{k_1+k_2}.\]
It is easy to see that $P^{n+1}$ satisfies the induction hypothesis.
\qed

\begin{rem}
The premises of Theorem \ref{thm:dihedral} are somewhat strong but necessary. First, the condition that the coefficients are real is necessary. To see this, consider \[P=3+\ii x-\ii x^{-1}+y.\]
In $\Z/3\Z\times \Z/2\Z$, we obtain
\[m_{\Z/3\Z\times \Z/2\Z}(3+\ii x-\ii x^{-1}+y) = \frac{\log 104}{6},\] 
however,
\[m_{D_3}(3+\ii x-\ii x^{-1}+y) =  \frac{\log 200}{6}.\]

It would be also nice to have a result that is more general for any $Q$ that is not necessarily reciprocal. Unfortunately,
this is not possible due to the different structures of the group rings.
Let
\[Q=\sum_{k=0}^{m-1} \alpha_k x^k  + \sum_{k=0}^{m-1} \beta_k y x^k.\]
We have, for $\Gamma_1 = \Z/m\Z\times \Z/2\Z$,
\[Q^*_1=\sum_{k=0}^{m-1} \overline{\alpha_k} x^{-k}  + \sum_{k=0}^{m-1} \overline{\beta_k} y x^{-k}.\]
For $\Gamma_2 = D_m$,
\[Q^*_2=\sum_{k=0}^{m-1} \overline{\alpha_k} x^{-k}  + \sum_{k=0}^{m-1} \overline{\beta_k} y x^k.\]

Clearly $QQ^*_1$ and $QQ^*_2$ are generally different:
\[QQ^*_1= \sum_{k_1,k_2} \left(\alpha_{k_1} \overline{\alpha_{k_2}} + \beta_{k_1} \overline{\beta_{k_2}}\right) x^{k_1-k_2} + \sum_{k_1,k_2} \left(\alpha_{k_1} \overline{\beta_{k_2}} + \beta_{k_1} \overline{\alpha_{k_2}}\right) y x^{k_1-k_2}\]

\[QQ^*_2= \sum_{k_1,k_2} \left(\alpha_{k_1} \overline{\alpha_{k_2}} + \overline{\beta_{k_1}}\beta_{k_2}\right) x^{k_1-k_2} + \sum_{k_1,k_2} \left(\beta_{k_1} \overline{\alpha_{k_2}} + \overline{\beta_{k_1}} \alpha_{k_2}\right) y x^{k_1-k_2}\]

As a counterexample for this case, let us look at $Q=x+2y$.

In $\Z/3\Z\times \Z/2\Z$, $QQ^*=5+2yx+2yx^{-1}$. Then
\[m_{\Z/3\Z\times \Z/2\Z}(x+2y) = \frac{\log 63}{6} .\] 

In $D_3$, $QQ^*=5+4yx^{-1}$ and
\[m_{D_3}(x+2y) = \frac{\log 3}{2}.\]

\end{rem}

Notice that the proof of Theorem \ref{thm:dihedral} is independent of $m$. Therefore, we can conclude the following.

\begin{coro}\label{coro:dihedral}
Let $P \in \R \left[\Z \times \Z/2\Z \right]$ be reciprocal. Then
\begin{equation}
m_{\Z\times \Z/2\Z}(P, \lambda) = m_{D_\infty}(P,\lambda),
\end{equation}
where $D_\infty = \langle \rho, \sigma \, | \, \sigma^2, \sigma \rho \sigma \rho \rangle $.
\end{coro}


\section{Dicyclic groups}

We consider the dicyclic groups with presentations
\[Dic_{m}=  \langle x, y \, | \, x^{2m},  y^2x^m, y^{-1}xyx \rangle.\]
The Mahler measure in these groups will be compared to the one in  $\Z/2m\Z \times \Z/2\Z = \langle x, y \, | \, x^{2m}, y^2, [x,y]\rangle$.

\begin{thm}
Let \[P=\sum_{k=0}^{2m} \alpha_k x^k + \sum_{k=0}^{2m} \beta_k y x^k\]
such that the coefficients $\alpha_k$ are real, and it is reciprocal in both $Dic_m$ and $\Z/2m\Z \times \Z/2\Z$. Then
\begin{equation}
m_{\Z/2m\Z \times \Z/2\Z} (P, \lambda)= m_{Dic_m} (P, \lambda).
\end{equation}
\end{thm}
\pf 
The proof is very similar to the one for Theorem \ref{thm:dihedral}. We just need to observe that the reciprocal condition in $Dic_m$ implies $\alpha_k=\alpha_{2m-k}$ and  $\beta_k=\overline{\beta_{m+k}}$, while in $\Z/2m\Z \times \Z/2\Z$ we have $\alpha_k=\alpha_{2m-k}$ and  $\beta_k=\overline{\beta_{2m-k}}$.
\qed

\section{Finite Approximations of the Mahler measure}
\label{sect:finite_approximation}

In this section we study a generalization of Theorem \ref{thm:approx} for families of groups that are not necessarily abelian.
Suppose that $\Gamma$ is a group and $\Gamma_m$ is a family of groups. Denote by
\[a_n^{(m)} = \left[P^n\right]_0 \qquad \mathrm{in} \quad \Gamma_m,\]
and let $a_n$ denote the corresponding number in $\Gamma$. Then

\begin{thm} \label{thm:generalapprox} Assume that there are $N_m$, depending on $m$ and $\{\Gamma_m,\Gamma\}$, with
$$\lim_{m \rightarrow \infty} N_m = \infty,$$
such that
\[ a_n^{(m)}= a_n \qquad \mathrm{for}\quad n < N_m.\]
Let $k$ be the $l_1$-norm of the polynomial $P$, i.e. the sum of the absolute values of the coefficients of the monomials in $P$.
Then, for $|\lambda| < \frac{1}{k}$,
\begin{equation}
\lim_{m \rightarrow \infty} m_{\Gamma_m}(P, \lambda) = m_{\Gamma}(P, \lambda).
\end{equation}
Moreover, the convergence is uniform in $\lambda$ for  $|\lambda|<\frac 1 k$.
\end{thm}

\pf Consider
\[\left| \sum_{n=1} \frac{a_n^{(m)} \lambda^n}{n} - \sum_{n=1} \frac{a_n \lambda^n}{n}\right|
\leq \left| \sum_{n=N_m} \frac{(a_n^{(m)}-a_n )\lambda^n}{n}\right|\]
for any $m$. The term in the right can be made arbitrarily small by 
choosing $m$ such that $N_m$ is large, since we know that the sums 
converge. 
Also, since $|\lambda|$ is bounded, the convergence is uniform. 

\qed

This Theorem provides an alternative way of seeing Theorem \ref{thm:approx} but it has other applications as well.
In the following corollary the $\Gamma_m$ are quotients of $\Gamma$:
  
\begin{coro}\label{everything} Let $P \in \Gamma$ reciprocal.
\begin{itemize}

\item For $\Gamma = D_\infty$, $\Gamma_m= D_m$,
\begin{equation} 
\lim_{m \rightarrow \infty} m_{D_m}(P, \lambda)= m_{D_\infty} (P, \lambda).
\end{equation}

\item For $\Gamma = Dic_\infty = \langle x, y \, | \, y^2, (yx)^2 \rangle $, $\Gamma_m = Dic_m$,
\begin{equation} 
\lim_{m \rightarrow \infty} m_{Dic_m}(P, \lambda)= m_{Dic_\infty} (P, \lambda).
\end{equation}

\item For $\Gamma = PSL_2(\Z) = \langle x, y \, |\, x^2, y^3  \rangle$, $\Gamma_m = \langle x, y \, |\, x^2, y^3, (xy)^m \rangle$,
\begin{equation} 
\lim_{m \rightarrow \infty} m_{\Gamma_m}(P, \lambda)= m_{PSL_2(\Z)} (P, \lambda).
\end{equation}

\item For $\Gamma = \Z * \Z = \langle x, y \rangle$, $\Gamma_m = \langle x, y \, |\, [x,y]^m \rangle$,
\begin{equation}
\lim_{m \rightarrow \infty} m_{\Gamma_m}(P, \lambda)= m_{\Z*\Z} (P, \lambda).
\end{equation}
Note that in this case $\Gamma_1=\Z \times \Z$.

\end{itemize}
\end{coro}

 

\section{Changing the base group}

In this section we fix an element in $\Z*\Z$, namely $P=x+x^{-1}+y+y^{-1}$ and we investigate $m_\Gamma(P,\lambda)$ for different choices of $\Gamma = \Z*\Z/N$.

\subsection{The free abelian group}
The computations of the first section show that for the case of $\Gamma = \Z\times\Z$, we get
\[ u_{\Z\times \Z}(P,\lambda) = \sum_{m=0}^\infty \binom{2m}{m}^2 \lambda^{2m} = F\left(\frac{1}{2},\frac{1}{2};1; 16\lambda^2\right) = \frac{2}{\pi} K(16 \lambda^2) ,\]
where $K$ is the complete elliptic function and $F$ the hypergeometric series. 

Recall that if 
\[u_{\Z\times \Z}(P,\lambda) = \sum_{n=0}^\infty a_n \lambda^n\]
is the generating function of the $\left[P^n\right]_0$, then
\[m_{\Z\times \Z}(P,\lambda) = -\sum_{n=1} \frac{a_n \lambda^n}{n}.\]

Therefore, for this case
\begin{equation}
 m_{\Z\times \Z}(P,\lambda) = -\sum_{m=1}^\infty \binom{2m}{m}^2 \frac{\lambda^{2m}}{2m}.
 \end{equation}

\subsection{Taking the limit over one variable}

Now assume that $\Gamma=\Z\times \Z/m\Z$ so that $x$ is free, $y^m=1$ and they commute. By using similar ideas to those in the proof of Theorem \ref{thm:approx}, we can see that
\begin{equation}
\lim_{n\rightarrow \infty} m_{\Z/n\Z\times\Z/m\Z}(P,\lambda)=m_{\Z\times\Z/m\Z}(P,\lambda).
\end{equation}

Therefore,
\begin{equation}
m_{\Z\times\Z/m\Z}(P,\lambda) = \frac{1}{m}  \sum_{k=0}^{m-1} m\left(1 - \lambda \left(x+x^{-1}+\xi_m^{k}+\xi_m^{-k}\right)\right),
\end{equation}
where the term in the right corresponds to the classical Mahler measure of one-variable polynomials.
The polynomial
\[ 1 - \lambda \left(x+x^{-1}+\xi_m^{k}+\xi_m^{-k}\right)\]
has two roots given by
\begin{equation} \label{eq:mixed}
\frac{\lambda^{-1}-2 \cos \frac{2\pi k}{m} \pm \sqrt{\lambda^{-2} -4 \cos \frac{2\pi k}{m} \lambda^{-1}- 4 \sin^2    \frac{2\pi k}{m}}}{2}.
\end{equation}

Since the product is $1$, one of the roots has always absolute values larger than $1$. Assume for simplicity that $\lambda>0$, then we need to consider the root with the "$+$" sign. 
Thus,
\begin{equation}
 m_{\Z\times\Z/m\Z}(P,\lambda) = \frac{1}{m}  \sum_{k=0}^{m-1}\log \left(\frac{1-2 \cos \frac{2\pi k}{m} \lambda + \sqrt{1 -4 \cos \frac{2\pi k}{m} \lambda- 4 \sin^2    \frac{2\pi k}{m}\lambda^2}}{2}\right).
 \end{equation}

For the case of $m=2$ the eigenvalues are $2$ and $-2$. Then we have, by equation (\ref{eq:mixed})
\[  m_{\Z\times\Z/2\Z}(P, \lambda) =-\log 2 + \frac{1}{2}  \left( \log \left(1-2 \lambda + \sqrt{1 -4 \lambda}\right) + \log \left(1+2 \lambda + \sqrt{1 +4\lambda}\right) \right).\]

In this case the coefficients can be also computed directly.
\[\left[\left(x+x^{-1}+y+y^{-1}\right)^n\right]_0 =  \sum_{j=0}^n \binom{n}{j}\left[(x+x^{-1})^j (y+y^{-1})^{n-j}\right]_0,\]
and because there is no relation between $x$ and $y$ besides the fact that they commute,
\[\left[(x+x^{-1}+y+y^{-1})^n\right]_0 = \sum_{j=0}^n \binom{n}{j}\left[(x+x^{-1})^j\right]_0 \left[(y+y^{-1})^{n-j}\right]_0.\]

Note that
\[\left[(x+x^{-1})^j\right]_0 = \left \{ \begin{array}{cc} \binom{j}{j/2} & j \quad \mbox{even}\\ 0 & j \quad \mbox{odd} \end{array} \right . \qquad  \left[(y+y^{-1})^j\right]_0 = \frac{2^j + (-2)^j}{2}.\]
Then
\[\left[\left(x+x^{-1}+y+y^{-1}\right)^{2l}\right]_0= \sum_{j=0}^{2l} \binom{2l}{2j} \binom{2j}{j} 2^{2l-2j}\]
\[ = \left[\left(x+x^{-1}+ 2 \right)^{2l}\right]_0 = \left[\left(x^2+x^{-2}+ 2 \right)^{2l}\right]_0 =
\left[\left(x+x^{-1}\right)^{4l}\right]_0 = \binom{4l}{2l}.\]
Then we obtain
\begin{equation}
m_{\Z\times\Z/2\Z}(P,\lambda) =  - \sum_{l=1}^\infty \binom{4l}{2l}\frac{\lambda^{2l}}{2l}.
\end{equation}

Let us remark that while we are computing a function of the form
\[\sum_{n=1}^\infty \sum_{j=0}^{n/2}\frac{1}{n} 
\binom{n}{2j} \binom{2j}{j} \beta^{2j}\lambda^n,\]
a similar function was studied in \cite{Dasbach:NaturalLog}:
\[\sum_{n=1}^\infty \sum_{j=0}^{n/2}\frac{1}{n} 
\binom{n}{j} \binom{2j}{j} \beta^{j}\lambda^n.\]

Finally, we remark, in the same vein as Theorem \ref{thm:approx}, that
\begin{equation}
\lim_{m \rightarrow \infty} m_{\Z \times \Z/m \Z}(P,\lambda) = m_{\Z \times \Z}(P,\lambda).
\end{equation}

\subsection{The free group}
We now assume that there is no relation between $x$ and $y$. Then $u_{\Z*\Z}(P,\lambda)$ corresponds to counting circuits (based in a distinguished point) in the $4$-regular tree (that is the Cayley graph of $\Z*\Z$). By a result of Bartholdi \cite{Bartholdi:Counting} $u$ is equal to
\[u_{\Z*\Z}(P,\lambda) = \frac{3}{1+2\sqrt{1-12\lambda^2}}.\]  

\subsection{$PSL_2(\Z)$}
A presentation for $PSL_2(\Z)$ is given by $\langle x,y \,|\,x^2, y^3\rangle$. Hence, it is isomorphic to $\Z/2\Z * \Z/3\Z$. By a result of Bartholdi \cite{Bartholdi:Counting},
 \[u_{PSL_2(\Z)}(x+y+y^{-1},\lambda) = \frac{(2-\lambda)\sqrt{1-2\lambda-5\lambda^2+6\lambda^3 +\lambda^4}-\lambda+\lambda^2+\lambda^3}{2(\lambda -1) (3 \lambda - 1) (2 \lambda + 1)}.\]
 

However, our $P$ is $2x+y+y^{-1}$. If we use Theorem 9.2 in \cite{Bartholdi:Counting} we obtain
\[u_{PSL_2(\Z)}(P,\lambda) = \frac{(2-\lambda)\sqrt{1-2\lambda-11\lambda^2+12\lambda^3 +4\lambda^4}-\lambda+\lambda^2-2\lambda^3}{2(\lambda-1)(3\lambda+1)(4\lambda -1)}.\]

\section{Arbitrary number of variables}
There is a general result in  \cite{Bartholdi:Counting} regarding the generating function of the circuits of a $d$-regular tree which is given by
\begin{equation}\label{eq:manyvariables}
g_d(\lambda)= \frac{2(d-1)}{d-2+d\sqrt{1-4(d-1)\lambda^2}}.
\end{equation}

In the case of 
\[P_{1,l}=x_1+x_1^{-1}+ \dots + x_l + x_l^{-1},\]
we can easily see that for the free group $\F_l$ we are  counting circuits in a $2l$-regular tree, therefore
\begin{equation}
u_{\F_l}(P_{1,l},\lambda) = g_{2l}(\lambda).
\end{equation}

We also consider the polynomial 
\[P_{2,l}=\left(1+x_1+\dots + x_{l-1}\right)\left(1+x^{-1}_1+\dots + x^{-1}_{l-1}\right).\]
Let us  add one more variable, so we obtain 
\[\left(x_0+x_1+\dots + x_{l-1}\right)\left(x_0^{-1}+x^{-1}_1+\dots + x^{-1}_{l-1}\right).\]
It is easy to see that the trace of the powers of $P_{2,l}$ is not affected by this change. But now it is clear that we are counting circuits in a $l$-regular tree, thus,
\begin{equation}
u_{\F_{l-1}}(P_{2,l},\lambda)= g_l(\lambda).
\end{equation}

 
In particular,
\begin{equation}
m_{\F_l}(P_{1,l}, \lambda) = m_{\F_{2l-1}}(P_{2,2l}, \lambda).
\end{equation}
 
Let us examine the abelian versions of these polynomials. By using elementary combinatorics we compute the trace of the powers of $P_{1,l}$ and $P_{2,l}$:
\[ \left[P_{1,l}^n\right]_0 = \sum_{a_1+\dots + a_l=n} \frac{(2n)!}{(a_1!)^2 \dots (a_l!)^2},\]
\[ \left[P_{2,l}^n\right]_0 = \sum_{a_1+\dots + a_l=n} \left(\frac{n!}{a_1! \dots a_l!}\right)^2.\]

For $l>2$ these expressions do not seem to simplify, in the sense that we are unable to find a closed formula that does not involve a summation. However, we make the following interesting observation:
\begin{equation}
\left[P_{1,l}^{2n}\right]_0 = \binom{2n}{n}\left[P_{2,l}^n\right]_0,
\end{equation}
where the trace in the left is taken over $\Z^{\times l}$ and the one in the right is taken over $\Z^{\times (l-1)}$.
 
In other words, there are relations among the Mahler measure of these two families of polynomials. But the relations depend on the base group. 
 
\bigskip 
 
{\bf Acknowledgements}: The authors would like to thank Neal Stoltzfus and Fernando Rodriguez-Villegas for helpful discussions. Part of this work was completed while ML was a member at the Institute for Advanced Study and during  visits to the Department of Mathematics at Louisiana State University. She is deeply grateful for their support and hospitality. 

\bibliography{linklit}
\bibliographystyle {amsalpha}

\end{document}